\documentclass[12pt]{article}
\usepackage{amsmath}
\usepackage{amsfonts}
\usepackage{amssymb}
\usepackage{theorem}

\title{ Sofic profiles of $S(\omega )$ and computability}
\author{A. Ivanov (Iwanow) 
}

\newtheorem{theorem}{Theorem}
\newtheorem{proposition}{Proposition}
\newtheorem{corollary}{Corollary}
\newtheorem{lemma}{Lemma}
\newtheorem{definition}{Definition}

\newtheorem{example}{Example}
\newtheorem{remark}{Remark}
\newtheorem{mtheorem}{Theorem} 
\newenvironment{proof}{\addvspace{8pt plus 2pt minus 2pt}\noindent\emph{Proof. }}
  { \begin{flushright}$\blacksquare$\par\addvspace{8pt plus 2pt minus
2pt}\end{flushright}}

\begin{document} 
\topmargin = 12pt
\textheight = 630pt 
\footskip = 39pt 

\maketitle

\begin{quote}
{\bf Abstract.}
We show that for every sofic chunk $E$ there is a bijective 
homomorphism $f: E_c \rightarrow E$, where $E_c$ is a chunk of 
the group of computable permutations of $\mathbb{N}$ so that the approximating morphisms of $E$ can be viewed as restrictions 
of permutations of $E_c$ to finite subsets of $\mathbb{N}$. 
Using this we study some relevant effectivity conditions 
associated with sofic chunks and their profiles.  
\end{quote}

\section{Introduction} 
\label{intro} 

We remind the reader that an abstract group $G$ is called 
{\bf sofic} if $G$ embeds into 
a metric ultraproduct of finite symmetric groups 
with the Hamming distance $d_{H}$, 
where the latter is defined as follows 
\[ 
d_H (g,h) = 1 - \frac{|\mathsf{Fix} (g^{-1}h )|}{n} \mbox{ , for } g,h \in S_n  
\]  
(see \cite{CL}, \cite{pestov}). 
It is an open question of M. Gromov 
whether every countable group is sofic. 

A more detailed approach to this notion was suggested by 
Y. de Cornulier in  \cite{cornulier2}. 
Firstly as in \cite{cornulier2} 
we call a finite subset $E$ of a group $G$ 
with $1$ a {\bf chunk} of $G$. 
In the definition below we assume that $\mathsf{inf} \emptyset = +\infty$. 

\begin{definition}
Let $E$ be a chunk of $G$ and let $n\in \mathbb{N}$, 
$0< \varepsilon < 1$. 
\begin{itemize} 
\item  
An $\varepsilon$-{\bf morphism} from $E$ to $(S_n, d_H)$ 
is a mapping $f: E \rightarrow S_n$ such that 
$f(1) = 1\in S_n$ and 
$d_H (f(xy), (f(x)f(y))) \le \varepsilon$ 
for all $x,y\in E$ with $xy\in E$. 
\item 
A mapping $f: E\rightarrow S_n$ is said to be 
$(1 -\varepsilon)$-{\bf expansive} if 
$d_H (f(x), f(y))\ge (1 - \varepsilon )$ 
whenever $x,y$ are distinct points of $E$.  
\item A $(1 - \varepsilon )$-expansive $\varepsilon$-morphism from $E$ is called a {\bf sofic approximation} of $E$ for $\varepsilon$.  
\item 
Define the {\bf sofic profile} of the chunk $E$ 
as  the function: 
\[
\mathsf{prof}_E (r) = \mathsf{inf }\{ n\in \mathbb{N} 
\mbox{: there exists a } (1 - r^{-1})\mbox{-expansive }
r^{-1}\mbox{-morphism} 
\] 
\[ 
E \rightarrow (S_n ,d_H ) \} \mbox{ , } r\in \mathbb{R} \mbox{ , } r>1.     
\] 
\item 
We say that the chunk $E$ is {\bf sofic} 
if for any $r>1$, $\mathsf{prof}_E (r) <\infty$.  
\item A group $G$ is {\bf sofic} 
if each chunk of $G$ is sofic. 
\end{itemize} 
\end{definition} 
 
As a result we have an ultraproduct-free definition of soficity  
(see also Remark 3.13 of \cite{cornulier2} where some comparison 
with the approach of Arzhantseva and Cherix \cite{AC} is given). 
We will consider chunks as partial groups. 
In particular we say that a map 
$f : E_1 \rightarrow E_2$ between two chunks is a {\bf homomorphism} 
if it preserves the multiplication and the unit. 
When $f$ is bijective we say that $E_1$ is a {\bf realization} of $E_2$. 
Note that in this case $f$ is not necessarily an isomorphism. 

Since permutation groups are involved into  the definition of soficity, 
they often arise in this context, see \cite{ES}, \cite{thomas}. 
In our case the idea is as follows. 
Note that every chunk can be realized in $S(\omega )$, the group of all 
permutations of the set of natural numbers. 
In order to study soficity we introduce a kind of approximating morphisms 
of such chunks. 

\begin{definition} \label{s-m} 
Let $E\subset S(\omega )$ be a chunk. 
Given $n\in \mathbb{N}$ a mapping $\sigma: E \rightarrow S_{n}$  
is called a {\bf supp-morphism} if  $\sigma (1) =1$ and 
\[ 
\forall \rho\in E \forall m< n (\rho (m) < n \rightarrow  
\rho (m) = \sigma (\rho )(m)) . 
\]  
\end{definition} 
The idea of supp-morphisms is to associate to each $\rho \in E$ 
a permutation from $S_n$  which can be considered as a restriction of 
$\rho$  to $\{ 0,\ldots ,n-1 \}$. 
It is obvious that any $\rho \in S(\omega )$ has a supp-morphic image in $S_n$. 

To recognize sofic chunks with morphisms supporting its soficity we will concentrate on chunks of $S(\omega )$ of some special form. 

\begin{definition} \label{g-chunk} 
Given  total function $g:\mathbb{N} \rightarrow \mathbb{N}$  
we will say that a chunk $E \subset S(\omega )$ is a $g$-{\bf chunk} if 
for every $\rho \in E$, every $n\in \mathbb{N}$ and $ m\le n$ 
we have $\rho (m) \le g(n)$. 
\end{definition} 
We also restrict the family of functions $g$ in this definition by the following condition. 
\begin{definition} \label{sl} 
We say that a total function $g:\mathbb{N} \rightarrow \mathbb{N}$ 
is {\bf slow} if $g(n+1) \geq g(n)> n$ for all $n \in \mathbb{N}$ 
and  
\[
\lim_{n\rightarrow \infty} \frac{n}{g(n)} = 1 .
\] 
\end{definition}
We now formulate our main theorem.

\begin{mtheorem} \label{mth}   
 A chunk $E$  is sofic if and only if  there is a computable 
slow function $g: \mathbb{N} \rightarrow \mathbb{N}$, 
a $g$-chunk  $E' \subset S(\omega )$ consisting of computable permutations and realizing $E$ under a bijective homomorphism $h: E' \rightarrow E$ 
such that there is a family of supp-morphisms 
$\{ \delta_n : n \in \mathbb{N}\}$  which provide sofic approximations of $E$ under the identification $h: E' \rightarrow E$. 

Moreover for a sofic chunk $E$ the profile of the property that 
$\delta_n$  are sofic approximations of $E$ is a computable function. 
\end{mtheorem} 
This theorem will be proved in Section 3 in a slightly stronger formulation. 
We will discuss there all notions appearing in it (for example see Remark \ref{prof_R}). 

Using Theorem \ref{mth} we study some relevant computability conditions 
associated with sofic chunks and their profiles. 
We mainly concentrate on computable complexity of some 
properties of sofic profiles like growth etc.  
This develops some issues of M. Cavaleri's investigations 
in \cite{cavaleri}, \cite{cavaleri1}, \cite{cavaleri2}.  
We also mention that since by Theorem \ref{mth} sofic chunks 
can be realized by computable permutations, there are  natural   
connections of our approach with previous work in computability theory  
concerning algebraic and computabily structure of the group of computable permutations, see \cite{morozov} and \cite{nies}. 

We finish this introduction by two definitions taken from \cite{cornulier2}.  
\begin{definition} \label{sim_pf} 
If $u,v: [1,\infty ] \rightarrow [0,\infty ]$ are non-decreasing functions, 
we write $u\preceq_{pf} v$ if there are positive real numbers 
$C,C',C''$ such that 
\[ 
(\forall r\ge 1 )(u(r) \le Cv(C' r) + C'' ),  
\] 
and we write $u \simeq_{pf} v$ if $u \preceq_{pf} v \preceq_{pf} u$. 
\end{definition} 

\begin{definition}  \label{gr-prof}
Let $G$ be an abstract group. 
The {\bf  sofic profile} of $G$ 
is the family of $\simeq_{pf}$-equivalence classes of the 
functions $\mathsf{prof}_E$ for all chunks $E$ of $G$.  
\end{definition} 
The latter definition explains the title of the paper. 
In fact below it will not be used explicitly. 
However to a large extent the paper concerns  
computable complexity of sofic profiles of 
computably enumerable groups which have embeddings into 
$S(\omega )$ with properties suggested by the formulation of Theorem \ref{mth} (see Section 4).

\section{Computability and computable profiles} 

We use standard material from the computability theory (see  
\cite{odifreddi} and \cite{sri}). 
In particular we use the universal computable function $\varphi (x,y )$, and the corresponding indexation of the family of all unary 
partial computable functions 
$\varphi_k (x) = \varphi (k,x)$, $k\in \mathbb{N}$. 
Then all $W_k = \mathsf{dom} (\varphi_k )$, $k\in \mathbb{N}$, 
form the standardly indexed family of all computably 
enumerable sets (\cite{sri}, Section 1).  
We assume that the reader knows the definition of the  
Arithmetical Hierarchy (\cite{odifreddi}, Section IV.1 and \cite{sri}, Section 4). 
The notation are taken from \cite{sri}. 
We remind the reader that a set $A$ is $\Sigma_n$-(resp. $\Pi_n$)-complete, 
if the problem of recognition of the members of $A$ belongs to 
$\Sigma_n$ (resp. $\Pi_n$) 
and any other problem of $\Sigma_n$ (resp. $\Pi_n$) is reducible to $A$. 

From now on we identify each finite set $F\subset \mathbb{N}$ with its  
G\" odel number.

\subsection{Computably enumerable groups} 

Let $G$ be a countable group generated by some $X\subseteq G$. 
The group $G$ is called {\bf recursively presented} (see Section 4.7 of \cite{LS}) if $X$ can be identified with $\mathbb{N}$ (or with some 
$\{ 0,\ldots , n\}$) so that $G$ has a computably enumerable set of relators in $X$.  
Below we give an equivalent definition, see Definition \ref{df1}. 
It is justified by a possibility of identification of the whole $G$ with $\mathbb{N}$. 
We follow the approach of \cite{khmi}. 
To simplify the exposition we always assume that $G$ is infinite. 

\begin{definition}\label{df0}
Let $G$ be a group and $\nu: \mathbb{N} \rightarrow G$ be a surjective function. 
We call the pair $(G,\nu)$ a \textbf{numbered group}.
The function $\nu$ is called a \textbf{numbering} of $G$.
If $g\in G$ and $\nu(n)=g$, then $n$ is called a number of $g$.
\end{definition}

\begin{definition}\label{df1}
A numbered group $(G,\nu)$ is \textbf{computably enumerable} if the set 
\[ 
\mathsf{MultT}:= \{(i,j,k)\mbox{ : } \nu(i) \cdot \nu(j)=\nu(k)\} 
\] 
is computably enumerable.
\end{definition}

\begin{definition}\label{comp_gr} 
A numbered group $(G,\nu)$ is \textbf{computable} 
if the set $\mathsf{MultT}$ as in Definition \ref{df1} is computable. 
\end{definition}
When $(G, \nu )$ is a computable numbered group and 
$\nu$ is a bijection then $\nu^{-1}$ is 
an isomorphism from $G$ to a group on $\mathbb{N}$ with computable 
multiplication.  
We say that this group is a \textbf{computable presentation} of $G$. 
  
Any finitely generated group with decidable word problem 
obviously has a computable presentation. 
This also holds in the case of the free group $\mathbb{F}_{\omega}$ 
with the free basis $\{ x_0 , \ldots , x_i , \ldots \}$. 
If we fix a computable presentation  of $\mathbb{F}_{\omega}$ 
and the corresponding numbering $\nu_F$, 
then for every recursively presented group $G= \langle X\rangle$ 
and a natural homomorphism $\rho : \mathbb{F}_{\omega} \rightarrow G$ 
(taking $\omega$ onto $X$) we obtain a numbering 
$\nu = \rho \circ \nu_F$ which satisfies Definition \ref{df1}. 

\begin{remark}\label{rp}
{\em Let $(G,\nu)$ be a computably enumerable group.} 

\begin{itemize} 
\item There exists a computable function $\Phi:\mathbb{N}\times \mathbb{N} \rightarrow \mathbb{N}$ such that for all $ x,y \in \mathbb{N}$ the equality $\nu(x)\nu(y) = \nu(\Phi(x,y))$ holds.
\item For every $x\in \mathbb{N}$ we can effectively find $y\in \mathbb{N}$ with $\nu(x)\nu(y)=1$.
\item The sets $\{n:\nu(n)=1\}$ and $\{(n_1,n_2): \nu(n_1)=\nu(n_2)\}$ are computably enumerable. 
\end{itemize} 
\end{remark}

\begin{remark} 
{\em If $(G, \nu )$ is a computable numbered group, then }
\begin{itemize}
\item the set $\{(n_1,n_2): \nu(n_1)=\nu(n_2)\}$ is computable; 
\item $G$ has a computable presentation (under another numbering if $\nu$ is not bijective). 
\end{itemize} 

\end{remark}

\subsection{Effective soficity}

Let $G$ be a countable group and $E$ be a chunk of $G$. 
Consider the profile $\mathsf{prof}_E$ as a function 
$\mathbb{N}\setminus \{ 0 \} \rightarrow \mathbb{N}$ 
(taking the restriction to $\mathbb{N}\setminus \{ 0\}$). 
Then we easily see the following lemma. 

\begin{lemma} \label{l_1}
If $E$ is  sofic then the function $\mathsf{prof}_E$ is computable. 
\end{lemma}

\begin{proof} 
The statement follows from the assumption that $\mathsf{prof}_E$ is total 
and the easy observation that set of pairs  
$\{ (n,m) : \mathsf{prof}_E (n) <m \}$ is computable.  
\end{proof}

\begin{itemize} 
\item Thus the $\simeq_{pf}$-class of total $\mathsf{prof}_E$ 
consists of subrecursive functions (i.e. bounded by a computable one). 
\end{itemize} 

{\bf The family of sofic chunks} can be viewed 
as the family 
\[ 
\Sigma\mbox{-}SC = \{ (E, k ) : E \mbox{ is a finite set with a partial binary operation and its unit,  } 
\] 
\[ 
\varphi_k 
\mbox{ is total non-decreasing} 
\mbox{ and } \forall n ( n>1 \rightarrow \mathsf{prof}_E (n)\le \varphi_k (n) )\} .  
\] 
Indeed, if $E$ is a finite set with a partial binary operation and its unit and $\mathsf{prof}_E$ is bounded above then $E$ is a chunk of 
some group (the corresponding metric ultraproduct). 
Moreover by Lemma \ref{l_1} if $\mathsf{prof}_E$ is bounded above 
by a computable function  then $\mathsf{prof}_E$ is computable.  
We note the following observations.  
\begin{itemize} 
\item The set of triples 
$\{ (E,n,m) : \mathsf{prof}_E (n) <m \}$ is decidable, 
where each $E$ is a finite set with a partial binary operation and its unit.  
\item If $\varphi_k$ is a fixed total and non-decreasing function  then  
\[ 
 \{ E : E \mbox{ is a finite set with a partial binary operation and its unit,  and }   
\] 
\[ 
\forall n (n>1 \rightarrow  \mathsf{prof}_E (n)\le \varphi_k (n) )\}   \in \Pi_1. 
\] 
\item Since the graph of the universal computable function $\varphi (x,y)$ 
is a computably enumerable set,  
\[ 
\mathsf{TotND} := \{ k : \varphi_k \mbox{ is total and non-decreasing } \} \in \Pi_2 . 
\] 
\end{itemize} 
These statements are easy exercises similar to ones given in Section 4.1 of \cite{sri}. 
Now  it is clear that $\Sigma$-$SC$ can be presented as 
a $\Pi_2$-set of natural numbers.   
Moreover  we have the following proposition. 

\begin{proposition} 
$\Sigma$-$SC$ is a $\Pi_2$-complete set.  
\end{proposition} 

\begin{proof} 
Note that $\mathsf{TotND}$ is reducible to $\Sigma\mbox{-}SC$. 
Indeed, let $ss(k)$ be the computable function such that for all $n$ 
the equality $\varphi_{ss(k)}(n) = \varphi_k (n) +2$ holds. 
Consider $\mathbb{Z}/2\mathbb{Z}$ as a a chunk. 
Then $\mathbf{prof}_{\mathbb{Z}/2\mathbb{Z}} (n) \le 2$. 
We see 
\[ 
k\in \mathsf{TotND} \Leftrightarrow (\mathbb{Z}/2\mathbb{Z}, ss(k) )\in \Sigma\mbox{-}SC. 
\] 
Thus to verify the statement of the proposition 
it remains to notice that $\mathsf{TotND}$ is $\Pi_2$-complete. 
Let $\mathsf{Inf} = \{ e \in \mathbb{N}$ : $W_e$ is infinite$\}$. 
This set is $\Pi_2$-complete, see \cite{sri} (Section 4.3.1). 
Thus it is enough to show that $\mathsf{Inf}$ is reducible to $\mathsf{TotND}$. 
To see the latter given $e\in \mathbb{N}$ we build a unary function $f$ by the following inductive algorithm. 
Assuming that $f$ is defined for all $i\le m$ we find the pair 
$(n,s)$ with the minimal Cantor number such that $f(m)<n$ and 
$\varphi_e (n)$ is computed within $s$ steps.   
Put $f(m+1) = n$. 
Then $e\in \mathsf{Inf}$ if and only if $f$ is total (and then it is increasing). 
The algorithm which given $e$ finds $t\in \mathbb{N}$ such that  
the Turing machine of $\varphi_t$ coincides with the one of $f$ 
gives the required reduction. 
\end{proof} 

In spite of this proposition sofic profiles of typical computable groups 
are uniformly computable in the sense of the following definition.  

\begin{definition} 
Let $(G, \nu )$ be a computably enumerable group.   
We say that $(G, \nu )$ is {\bf effectively sofic}   if there is a uniform algorithm which for every finite $D\subset \mathbb{N}$ 
and every $n>1$ finds the value $m = \mathsf{prof}_E (n)$ for $E=\nu (D) \cup\{ 1\}$. 
\end{definition} 
It is worth noting that when $(G, \nu )$ has an algorithm as in this 
definition then producing $m = \mathsf{prof}_E (n)$ we can also 
algorithmically find the corresponding $(1 - n^{-1})$-expansive $n^{-1}$-morphism into $S_m$ (finding an appropriate one in a finite list of morphisms). 

In fact effective soficity was introduced by M. Cavaleri in \cite{cavaleri} 
in the case of finitely generated groups. 
We now formulate a version of Theorem 3.3.1 from \cite{cavaleri}.  

\begin{proposition} 
A computably enumerable group $(G, \nu )$ is effectively sofic if and only if  
$(G, \nu )$ is sofic and computable.  
\end{proposition}  

\begin{proof} 
It is clear that every effectively sofic group is sofic. 
To verify the equality $\nu (i) \cdot \nu (j) = \nu (k)$ 
let us apply effective soficity to $D=\{ i,j,k\}$ and $n=3$.  
Computing in the corresponding $S_m$ the distance between 
$\nu (i) \cdot \nu (j)$ and $\nu (k)$ we check if it is $< \frac{1}{3}$ or 
$\ge \frac{2}{3}$. 

If $(G, \nu )$ is sofic and computable, $n\in \mathbb{N}\setminus \{ 0,1\}$ and $D$ is a finite subset of $\mathbb{N}$ we start the procedure verifying for natural numbers $1,\ldots ,m, \ldots$ if there is a 
$(1 - n^{-1})$-expansive $n^{-1}$-morphism into $S_m$. 
This gives an algotithm for effective soficity. 
\end{proof}

\section{Function growth and the corresponding subgroups of $S(\omega )$}

\subsection{Growths and permutation groups} 

We will denote $\mathbb{N}_{\infty}=\mathbb{N} \cup \{ \infty \}$ and assume that $n < \infty$ for all $n \in \mathbb{N}$. 
\begin{definition} 
Let $\Omega_{\infty}$ be the semigroup of all total functions 
$g:\mathbb{N}_{\infty} \to \mathbb{N}_{\infty}$ 
such that $g(\infty)=\infty$ and 
$g(n+1) \geq g(n)> n$ for all $n \in \mathbb{N}$ . 
\end{definition} 
Clearly, the semigroup $\Omega_{\infty}$ has no identity element. 
Any total function $f: \mathbb{N} \rightarrow \mathbb{N}$ can be considered as a function $\mathbb{N}_{\infty} \to \mathbb{N}_{\infty}$ with 
$f(\infty) = \infty$. 
If the latter one is in $\Omega_{\infty}$ we say that $f\in \Omega_{\infty}$. 

We order $\Omega_{\infty}$ as follows: 
\[ 
f \prec g \Leftrightarrow \mbox{ there exists } n_0 \in \mathbb{N} 
\mbox{ such that } \forall n>n_0  (f(n) < g(n)) , 
\] 
and define 
\[ 
f \sim g \Leftrightarrow (\exists k \in \mathbb{N}) (f \prec g^k \wedge g \prec f^k ), 
\]
where $f^k$ is the $k$-th power of $f$ with respect to composition.   
The equivalence classes of $\sim$ will be called {\bf growths}.  

We say that  $g\in \Omega_{\infty}$ is {\bf slow} if it is 
slow on $\mathbb{N}$ (see Definition \ref{sl}). 
\begin{lemma} \label{grow}  
(a) Each growth $[f]_{\sim}$ is a subsemigroup of $\Omega_{\infty}$.\\ 
(b) If a growth $[g]_{\sim}$ contains a slow function 
then all functions of the growth are slow.
\end{lemma}  

\begin{proof} 
(a) This is easy and already is proved in \cite{hol}.   \\   
(b) This statement follows from the observation that given functions 
$g,h\in \Omega_{\infty}$ satisfying  
\[
\lim_{n\rightarrow \infty} \frac{n}{g(n)} = 1 \mbox{ and } 
\lim_{n\rightarrow \infty} \frac{n}{h(n)} = 1 .  
\] 
we also have: 
\[
\lim_{n\rightarrow \infty} \frac{n}{(gh)(n)} = 
\lim_{n\rightarrow \infty} \frac{n}{h(n)}\cdot \frac{h(n)}{gh(n)} = 1 .  
\] 
\end{proof} 
By statement (b) of the lemma it is natural to use the term 
{\bf slow growth} for any growth $\alpha$ containing a slow function. 
Note that the growth represented by $f(n) = n+1$ is slow.  

Slow growths will usually appear in the following context. 

\begin{definition}   
A permutation $\rho \in S(\omega )$ is {\bf  bounded} by $g \in \Omega_{\infty}$ if   
for every $n\in \mathbb{N}$ and $ m\le n$ we have $\rho (m) \le g(n)$.  
\end{definition}  
It is worth noting that according Definition \ref{g-chunk} 
a chunk $E \subset S(\omega )$ is a $g$-{\bf chunk} 
if and only if all elements of $E$ are bounded by $g$. 
It is obvious that: 
\begin{itemize} 
\item if $\rho_1$ is bounded by $g_1$ and $\rho_2$ is bounded by $g_2$ then 
$\rho_1 \cdot \rho_2$ is bounded by $g_1\circ g_2$.  
\end{itemize} 

\begin{definition} 
Let $\alpha = [f]_{\sim}$ be a growth.  
We define 
\[ 
S_{\alpha}(\omega ) = \{ \rho \in S(\omega ): \rho \mbox{ and } \rho^{-1} 
\mbox{ are bounded by some }g \in \alpha\}.
\]   
\end{definition} 
The  chunks of $S_{\alpha} (\omega)$ and their corresponding profiles are 
the main objects of our further investigations.
Let us mention the following remarks. 

\begin{itemize} 
\item   $S_{\alpha}(\omega )$ is a subgroup of $S(\omega )$ containing $SF(\omega )$, 
the subgroup of all finitary permutations.   
 
\item If $\alpha$ is represented by $f=\infty$, then  $S_{\alpha}(\omega ) = S(\omega )$. 
\end{itemize}

The following statement is a corollary of Lemma \ref{grow}.

\begin{corollary} \label{profiles} 
Let $\alpha$ be a slow growth. 
Then for any chunk $E\subseteq S_{\alpha}(\omega)$ there is a slow 
$g\in \alpha$ such that all elements of $E$ are bounded by $g$, i.e. $E$ is a $g$-chunk. 
\end{corollary}

\subsection{Profiles of $g$-chunks}

\begin{definition} \label{prof_g} 
Define the {\bf profile} of  $g\in \Omega_{\infty}$ as 
\[
\mathsf{prof}_g (r) = g(\mathsf{inf} \{ m\in \mathbb{N} : \frac{m}{g(m)} > 1 - r^{-1} \} )\mbox{ , } r> 1. 
\] 
\end{definition} 
\begin{example} \label{Ex1} 
{\em If $g(n) = n + [\sqrt{n}]$, then 
$\mathsf{prof}_g (r) = ((r-1)^2 +1) + r-1 = (r-1)^2 +r$. } 
\end{example}

\begin{remark} \label{prof_R} 
{\em Let $R \subseteq \mathbb{N} \times \mathbb{R}_{>1}$ be a binary relation. 
Then one can associate to $R$ the {\bf profile} of the property $R$ 
defined to be the following function $\mathbb{R}_{>1} \rightarrow \mathbb{N}$: 
\[ 
\mathsf{prof}_R (r) = \mathsf{inf} \{ n \in \mathbb{N} : (n,r) \in R\} 
\mbox{ , } r\in \mathbb{R} \mbox{ , } r>1 .  
\] 
Let us observe that Definition \ref{prof_g} can be reformulated in this way. 
Indeed,  
\[
\mathsf{prof}_g (r) = \mathsf{inf} \{ n\in \mathbb{N} : (\exists m \in \mathbb{N})( g(m) \le n )\wedge (\frac{n-m}{n} < r^{-1} )\} \mbox{ , } r> 1   
\]  
(note that the value of the right part of this equation 
is of the form $g(m)$). 
We see that this is $\mathsf{prof}_R$ where }  
$R (n,r) = (\exists m \in \mathbb{N})( g(m) \le n )\wedge (\frac{n-m}{n} < r^{-1} )$. 
\end{remark} 

It is easy to see that:  
\begin{itemize} 
\item when $f\prec g$, then $\mathsf{prof}_f \prec \mathsf{prof}_g$.  
\item if a function from $\Omega_{\infty}$ is slow then  
its profile is a total function. 
\end{itemize} 

The following proposition compares  
$\mathsf{Prof}_E$ and $\mathsf{Prof}_g$, where $E$ is a $g$-chunk. 
In the formulation we use supp-morphisms defined in Definition \ref{s-m} and the convention of Remark \ref{prof_R}. 
 
\begin{proposition} \label{vs}
Let $g: \mathbb{N} \rightarrow \mathbb{N}$ have total $\mathsf{prof}_g$ 
and  let $E \subset S(\omega )$ be a $g$-chunk. 
For every $n\in \mathbb{N}$ choose a supp-morphism $\sigma_n: E \rightarrow S_n$.  

Then the following statements hold. 
\begin{enumerate} 
\item  For any $\varepsilon$ there is $n$ 
such that $\sigma_n$ is an $\varepsilon$-morphism. 
\item  The profile of the property that the mapping 
$\sigma_n$  
is an $r^{-1}$-morphism 
is bounded above by $\mathsf{prof}_g (2r)$. 
\item  If for any $r$ and $n\in \mathbb{N}$ 
with $\mathsf{prof}_g (r) = n$ and any  
distinct $\rho_1 , \rho_2\in E$ we have  
$$ 
g(n-|\mathsf{Fix}(\sigma_{n} (\rho_1 )(\sigma_{n} (\rho_2))^{-1})|) \ge n,  
$$  
then the sofic profile of $E$ is bounded above by $\mathsf{prof}_g (2r)$.   
\end{enumerate} 
\end{proposition} 

\begin{proof} 
 By Definition \ref{s-m} 
\[ 
\forall \rho\in E \forall m\le n (\rho (m) \le n \rightarrow  
\rho (m) = \sigma_n (\rho )(m)) . 
\]   
Since $E$ is a $g$-chunk  it is easy to see that for any $\rho \in E$ 
the $\sigma_{n}$-image of $\rho$ can differ from 
the partial map which is the restriction $\rho |_{n}$, 
only inside $\{ m, ...,n \}$ where $m$ is maximal with $g(m) \le n$. 
Thus $\sigma_{n}(\rho )\sigma_{n}(\rho' )$
can differ from $(\rho \rho' )|_{n}$ only inside 
$g^{-1} (\{ m, \ldots , n \}) \cup \{ m, \ldots , n \}$. 
So by the definition of the Hamming metric, 
\[
d_H (\sigma_{n}(\rho )\sigma_{n}(\rho' ), \sigma_{n}(\rho\rho' )) \le \frac{2(n - m)}{n}. 
\]
Now statement 1 of the proposition follows from the definition of $\mathsf{prof}_g$. 
 
Statement 2 follows from the proof of statement 1. 

To see statement 3 note that the inequality 
\[
n - |\mathsf{Fix} (\sigma_{n} (\rho_1 )\sigma_{n} (\rho_2 )^{-1})| \ge m
\] 
implies that 
$d_H (\sigma_{n} (\rho_1 ),\sigma_{n} (\rho_2)) \ge 1 - \frac{n-m}{n}$. 
If $\frac{n - m}{n} \le (2r)^{-1}$, 
then  
\[ 
d_H(\sigma_{n} (\rho_1 ),\sigma_{n} (\rho_2 )) \ge 1 - (2r)^{-1}. 
\]  
The rest follows from statement 2.  
\end{proof}

\begin{corollary} \label{profiles} 
Let $\alpha$ be a slow growth. 
\begin{enumerate} 
\item 
Then for any chunk $E\subseteq S_{\alpha}(\omega)$ there is  $g\in \alpha$ 
such that $E$ is a $g$-chunk and for any choice of 
supp-morphisms $\sigma_n: E \rightarrow S_n$, $n\in \mathbb{N}\setminus \{ 0,1\}$,  the profile of the property that the mapping 
$\sigma_n$ is an $r^{-1}$-morphism  is bounded above by  $\mathsf{prof}_g (2r)$. 
\item Moreover given the multiplication table of $E$ the profile of the property that the mapping $\sigma_n$ is an $r^{-1}$-morphism  is a computable function on $\mathbb{N}$. 
\end{enumerate}
\end{corollary} 

\begin{proof} 
The first statement of the corollary easily follows from 
the definition of $S_{\alpha} (\omega)$, 
Lemma \ref{grow} and Proposition \ref{vs}.  
The second statement of the corollary follows from the fact that the corresponding profile is total 
and the easy observation that set of pairs  
$\{ (n,m) : \sigma_n$ is an $m^{-1}$-morphism $\}$ is computable.  
\end{proof}

It is worth noting that the profile of a slow computable function is computable. 
Thus if $\alpha$ above is represented by a computable function, 
then the corresponding  $\mathsf{prof}_g (2r)$ is computable on 
$\mathbb{N}$.

\begin{remark}  
{\em There is no reason why $g$ is slow under the assumption that 
$\mathsf{prof}_g$ is total. 
This would be the case if we replace $\mathsf{prof}_g$ 
by the profile, say $\mathsf{prof}^*_g$, of the following relation 
(see Remark \ref{prof_R}): 
\[ 
R^* (n,r) := (\forall k \ge n)(\exists m \in \mathbb{N})( g(m) \le k )\wedge (\frac{k-m}{k} < r^{-1} ). 
\] 
However note that such a replacement makes Proposition \ref{vs} weaker. 
Additionally there is no reason why a slow computable function $g$ has computable profile of the relation $R^* (n,r)$. 
We only mention the following curious observation concerning $R^*$. }
\end{remark} 

\begin{proposition} 
 If $g$ and $h$ are slow and have the same growth, then 
$\mathsf{prof}^*_g$ and $\mathsf{prof}^*_h$ are $\simeq_{pf}$-equivalent.  
\end{proposition} 

\begin{proof}  
To prove the proposition it is enough to show that 
$\mathsf{prof}^*_{g}$ and $\mathsf{prof}^*_{g^2}$ are 
$\simeq_{pf}$-equivalent. 
It is clear that for any $C>1$, $g(n) \prec C\cdot n$. 
Choosing sufficiently large $C$ we may additionally assume that 
$g(n) < C\cdot n$ for all $n>0$. 
If for some $n$ we have the property that for any $l\ge n$ 
there is a number $m$ with $g(m) \le l$ and 
$|\frac{m}{l} - 1|<\varepsilon$, then for each $k\ge C\cdot n$ 
we find two numbers $m_1 \ge [\frac{k}{C}]$ and $m_2$ such that 
$g(m_1 ) \le k$, $g(m_2 ) \le m_1$, 
\[ 
|\frac{m_1}{ k} - 1|<\varepsilon 
\mbox{ and } 
|\frac{m_2}{m_1} - 1| < \varepsilon.  
\] 
Thus 
\[  
|\frac{m_1}{ k} \cdot \frac{m_2}{m_1} - \frac{m_1}{ k}| < \varepsilon  
\] 
\[  
\mbox{ and } |\frac{m_2}{k} - 1 |< 2\varepsilon. 
\]  
We conclude that 
\[  
(\forall r \ge 1)(\mathsf{prof}^*_{g^2}(r) \le C \mathsf{prof}^*_{g}(\frac{1}{2} r) ) .  
\]  
\end{proof} 

\subsection{Realization} 

The following theorem is the main statement of this section. 
It justifies $g$-chunks of slow functions $g$. 

\begin{theorem} \label{realization} 
A chunk $E$  is sofic if and only if  there is a slow $g\in \Omega_{\infty}$ and
a $g$-chunk  $E' \subset S(\omega )$ realizing $E$ under a bijective 
homomorphism $h: E' \rightarrow E$ 
such that there is a family of supp-morphisms 
$\{ \delta_n : n\in \mathbb{N}\}$  which  provide sofic approximations of $E$ under the identification $h: E' \rightarrow E$. 

Moreover for a sofic chunk $E$ the function  $g$, 
the $g$-chunk $E'$ as above and the profile of the property that 
$\delta_n$, $n\in \mathbb{N}$, are sofic approximations of $E$ 
can be realized by computable functions. 
\end{theorem} 

\begin{proof} 
For every $n>1$ fix $m_n = \mathsf{prof}_E (n)$ and an appropriate  
$(1 - \frac{1}{n})$-expansive $\frac{1}{n}$-morphism 
$\sigma_n : E \rightarrow (S_{m_n}, d_H )$. 
Let $f$ be a function $\mathbb{N} \rightarrow \mathbb{N}$. 
The direct sum $\bigoplus_{1 <i\le n} S^{f(i)}_{m_i}$ is naturally 
considered as a permutation group on the set of natural numbers 
less than $\sum_{1< i\le n} f(i)m_i$.  
We associate to every $e\in E$ and $n\in \mathbb{N} \setminus \{ 0,1\}$ the permutation 
\[ 
\sigma_{\oplus,n} (e) = \bigoplus_{1 <i\le n} \sigma^{f(i)}_{i}(e) \in \bigoplus_{1 <i\le n} S^{f(i)}_{m_i}. 
\]   
Let $\sigma (e)$ be the limit of these elements 
in the natural direct limit of all $\bigoplus_{1 <i\le n} S^{f(i)}_{m_i}$, 
$n \rightarrow \infty$.  
Then $\sigma (e) \in S(\omega )$, the permutation $\sigma_{\oplus,n} (e)$ is the restriction of $\sigma (e)$ to the initial segment of $\sum_{1 <i \le n} f(i) m_i$ elements of $\mathbb{N}$ and the distance between $\sigma_{\oplus,n} (e)$ and $\sigma_{\oplus,n} (e')$ 
in $S_{\sum_{1 <i\le n} f(i)m_i}$ is equal to 
\[ 
\frac{\sum_{1 <i \le n} f(i) (m_i - |\mathsf{Fix} ( \sigma_{i} (e) \sigma_{i} (e')^{-1})|)}{\sum_{1 <i \le n} f(i) m_i} .
\] 
To satisfy the conditions of the theorem, i.e. to realize $\delta_n$ as 
$\sigma_{\oplus ,n}$, we need an appropriate function $f$. 
We define it by induction. 
Given $f(2), \ldots ,f(n-1 )$ one can choose $f(n)$ so that  $\sigma_{\oplus,n}$ 
is an $(1 - \frac{1}{n-1})$-expansive $\frac{1}{n-1}$-morphism into 
$S_{\sum_{1 <i\le n} f(i)m_i}$. 

We define the function $g$ as follows. 
Let $g(j) = j +  m_2 + \ldots + m_n$ for each $j$ with 
$\sum_{1 <i < n} f(i) m_i \le j < \sum_{1 <i \le n} f(i) m_i$ 
(where the empty sum is $0$).  
Then each $\sigma (e)$ with $e\in E$ is bounded by $g$. 
Thus $E' = \{ \sigma (e) : e\in E\}$ is a $g$-chunk. 

Given $f(2), \ldots ,f(n-1 )$ one can correct $f(n)$ so that 
\[ 
\frac{\sum_{1 <i < n} m_i }{(\sum_{1 <i \le n} f(i) m_i) -1 + \sum_{1 <i < n} m_i } <\frac{1}{n} . 
\] 
Since this exactly means 
\[ 
1 - \frac{(\sum_{1 <i \le n} f(i) m_i )-1 }{g((\sum_{1 <i \le n} f(i) m_i )-1)} <\frac{1}{n} ,   
\] 
we conclude that taking $f$ sufficiently fast 
we guarantee that $g$ is slow and the maps $\sigma_{\oplus,n}$ 
form a family of supp-morphisms which are 
$(1 - \frac{1}{n-1})$-expansive $\frac{1}{n-1}$-morphisms into 
$S_{\sum_{1 <i\le n} f(i)m_i }$. 

Since $E$ is sofic, by Lemma \ref{l_1} the sequences $m_2 , m_3 ,\ldots ,m_n ,\ldots$ is computable. 
In particular we can define an algorithm of computation of $f$ 
so that the inequations guaranteeing slowness of $g$ 
and the property that $\sigma_{\oplus,n}$ 
is an $(1 - \frac{1}{n-1})$-expansive $\frac{1}{n-1}$-morphism into 
$S_{\sum_{1 <i\le n} f(i)m_i}$, are satisfied (verifying natural numbers in turn).
Computability of permutations $\sigma(e)$, $e\in E$, follows from 
computablity of the sequence $\sigma_2 , \sigma_3 , \ldots ,\sigma_n ,\ldots$ 
in the beginning of the proof. 
\end{proof}  
 
Theorem \ref{realization} can be interpreted as an evidence of some connection between soficity and computability.  
In some sense it says that soficity of a chunk is equivalent 
to possibility of its computable realization in $S(\omega )$. 
Thus it can be viewed as an effective version of the observation that each chunk can be realized in $S(\omega )$  
(the latter is an obvious consequence of the fact that 
each countable group  embeds into $S(\omega )$).  
 
\begin{remark} 
{\em Consider finitely generated groups with the word problem in $\Pi_1$.  
It is easy to see that any finitely generated group of computable permutations belongs to this class. 
Answering a question of Higman,  Morozov has proved in \cite{morozov} that there is a 2-generated group with the word problem in $\Pi_1$
which is not embeddable into the group of computable permutations. 
This is slightly opposite to Theorem \ref{realization}. 
It is natural to ask if the Morozov's example is sofic or if there is a sofic group which witnesses Morozov's theorem. 
Analysing the construction of \cite{morozov} one easily notes that under the conjecture that HNN-extensions preserve soficity, Morozov's example is sofic. 
The former conjecture is an open question.  }
\end{remark} 

\begin{remark} 
{\em It is worth noting that Theorem \ref{realization} does not state that a sofic chunk is contained in some $S_{\alpha} (\omega )$ with a slow $\alpha$. 
In particular we do not know if a sofic group must be locally embeddable into the class of all $S_{\alpha} (\omega)$ with slow $\alpha$. 
The author is grateful to the referee for this question. }  
\end{remark} 

The following example is slightly connected with 
Problem 3.18 from \cite{cornulier2} which asks if there is a group 
for which the sofic profile is unbounded and not linear.  
Indeed assuming below that for the slow function $g$ the profile 
$\mathsf{prof}_g (2r)$ is non-linear, we can choose a family of supp-morphisms  $\sigma_n$, $n\in \mathbb{N}$, so that the $g$-chunk $E$ below has the property  that  linearity of its sofic profile cannot be demonstrated by these $\sigma_n$. 
Note that Example \ref{Ex1} gives a concrete slow $g$ with non-linear profile.   

\begin{example}
For any slow growth $\alpha$ there is function $g\in \alpha$, a  
$g$-chunk $E\subseteq S_{\alpha}(\omega)$ and a family of supp-morphisms 
$\sigma_n: E \rightarrow S_n$, $n\in \mathbb{N}$,  such that the profile of the property that $\sigma_n$  is an $r^{-1}$-morphism  
is equivalent to $\mathsf{prof}_g (2r)$.  \\ 
{\em We build a 3-element chunk as follows. 
Let $h$ be a disjoint union of $3$-element cycles: 
$$ 
\mbox{ for each } k\in \mathbb{N} 
\mbox{ let } h(3k) = 3k+1 \mbox{ , } 
h(3k+1) = 3k+2 \mbox{ , and } h(3k+2) = 3k. 
$$ 
Let $E$ consist of $h, h^2$ and $\mathsf{id} = h^3$. 
Let $g\in \alpha$ be sufficiently fast, for example 
$g(n) >30 +n$. 
For $n= 3m+ i$ with $i\in \{ 0, 1 \}$ let $\sigma_n (h)$ 
coincide with $h$ for all $l <3m$ and $\sigma_n (h^2 )$ 
coincide with $h^2$ for all $l$ with $g(l) \le n$. 
We additionally assume that for all 
$l$ with $g(l-2) > n$ the permutation 
$\sigma_n (h^2 )$ coincides with $\mathsf{id}$.  
Then assuming that $m$ is the maximal number with $g(m)\le n$ we have 
$|m - |\mathsf{Fix} ((\sigma_n (h))^{-2} \sigma_n (h^{2}))|| \le 5$. 
This implies that the profile of the property that $\sigma_n$  is an $r^{-1}$-morphism  is equivalent to $\mathsf{prof}_g (2r)$. } 
\end{example}

\subsection{$g$-Chunks of computable functions} 

In this subsection we show that Theorem \ref{realization} 
provides a parametrization of sofic chunks which is 
different from the one given in Section 2.2,  
but which is still $\Pi_2$-complete. 

\begin{theorem} \label{3_2} 
The family of all triples $(D,(\circ ,e),n)$ with conditions given below 
is $\Pi_2$-complete:  
\begin{itemize}
\item $n \in \mathbb{N}$ such that $\varphi_n \in \Omega_{\infty}$ and the profile of $\varphi_n$ is total,
\item $D$ is a finite subset of $\mathbb{N}$, $\circ$ is a partial binary operation on $D$ where $e$ is neutral, and the set 
$\{ \varphi_k  : k\in D\}$ forms a $\varphi_n$-chunk of $S(\omega )$ which realizes $(D, \circ ,e)$ with respect to the map $\varphi_k \rightarrow k$,  and 
\item the chunk $(D, \circ ,e)$ is sofic under approximations by 
supp-morphisms. 
\end{itemize} 
\end{theorem}

\begin{proof} 
Let $\mathcal{R}$ be the the following family 
\begin{itemize} 
\item all pairs $(D, (\circ ,e))$, where $D$ is a finite subset $\mathbb{N}$, $\circ$ is a partial binary operation on $D$ where $e$ is neutral such that the set $E= \{ \varphi_k : k\in D\}$ 
forms a chunk of $S(\omega )$ with $\varphi_e = \mathsf{id}$, which realizes the chunk $(D, \circ ,e)$ under the map $\varphi_k \rightarrow k$.  
\end{itemize} 
To see that $\mathcal{R} \in \Pi_2$ we associate to each chunk 
$(D, \circ ,e)$ the set $\Theta(D, \circ ,e)$ consisting of all conditions of the form 
\[ 
k \not= k' \mbox{ , } k \circ k' \not= k'' , 
\]
which are satisfied in $(D, \circ , e)$. 
It is easy to see that the following set is computable: 
\begin{itemize} 
\item  all quadruples $((D,(\circ , e), f,s)$ where $f$ is a function 
$\Theta(D, \circ ,e)\rightarrow \mathbb{N}$ such that all conditions 
\[ 
\varphi_k (f(\theta ))\not= \varphi_{k'} (f(\theta )) 
\mbox{ , } \varphi_k \circ \varphi_{k'} (f(\theta ))\not= \varphi_{k''} (f(\theta ))  
\] 
are satisfied after $s$ steps of computation of terms appearing in them,  where $\theta$'s denote the inequalities from $\Theta(D, \circ ,e)$ in which the corresponding $f(\theta )$'s  are substituted.
\end{itemize} 
Taking the projection to the components $(D, (\circ ,e))$ we obtain a computably enumerable set containing $\mathcal{R}$. 
Indeed, when $(D, (\circ ,e))\in \mathcal{R}$, the map $\varphi_k \rightarrow k$ is a 1-1-homomorphism of chunks, i.e. every inequality from $\Theta(D, \circ ,e)$ is satisfied by $\varphi_k$'s at some natural number (and this defines the corresponding value of $f$ and the step $s$ of computation).  

To finish the proof that $\mathcal{R} \in \Pi_2$ it remains to show 
the set of all pairs $(D,e)$ such that the set $\{ \varphi_k  : k\in D\}$ consists of permutations and $\varphi_e = \mathsf{id}$, belongs to $\Pi_2$.
This follows from computable enumerability of the graph of the universal 
computable function $\varphi (x,y)$ and appropriate $\forall$-quantifying.  

Let $\mathcal{R}^*$ be the following family
\begin{itemize}   
\item all triples $(D,(\circ ,e), m )$, where  
$(D, (\circ ,e))\in \mathcal{R}$ and the chunk $\{ \varphi_k : k\in D\}$ 
forms a $g$-chunk where $g = \varphi_m \in \Omega_{\infty}$ and the profile of $g$ is total. 
\end{itemize} 
To see that $\mathcal{R}^* \in \Pi_2$ note that the set 
\[
\{ m :\varphi_m \in \Omega_{\infty} \mbox{ and } \forall n \exists n_1  (n \cdot n_1 > (n-1) \varphi_m (n_1 )) \}  
\] 
belongs to $\Pi_2$. 
This set consists of indicies of $\varphi_m \in \Omega_{\infty}$ with total profiles. 
Moreover, 
\[ 
\{ (k,m): \varphi_m\mbox{ is total and } \varphi_k\mbox{ and } \varphi^{-1}_k\mbox{ are bounded by }\varphi_m \} \in \Pi_2. 
\] 
Since $\mathcal{R} \in \Pi_2$ we easily deduce $\mathcal{R}^* \in \Pi_2$. 

We now discuss approximations by supp-morphisms. 
Note that given $\sigma: D \rightarrow S_n$ and chunk 
$E= \{ \varphi_k : k\in D\}$ of $S(\omega )$ to express that 
$\sigma$ defines an $(1-\frac{1}{l} )$-expansive 
$\frac{1}{l}$-supp-morphism of $E$ we only need existential 
quantifiers: for values $\varphi_k (s)$ with $s<n$ and for values  
supporting conditions of the form 
$\varphi_i \circ \varphi_j \not= \varphi_k$, where $i,j,k \in D$. 
This naturally provides a computably enumerable set of tuples 
$(D,  (\circ ,e),l, n )$. 
When a tuple from this set can be extended by a component $m$ with 
$(D, (\circ ,e), m) \in \mathcal{R}^*$ then  there exists a supp-morphism $\sigma_n$ to $S_n$ which is a $(1 - \frac{1}{l})$-expansive $\frac{1}{l}$-morphism. 
Using this and the fact that $\mathcal{R}^* \in \Pi_2$ we obtain 
that the following family belongs to the class $\Pi_2$:   
\begin{itemize}
\item The family of all triples $(D,(\circ ,e), m)$, 
where $(D, (\circ ,e), m)\in \mathcal{R}^*$ 
and the chunk $(D, \circ ,e)$ is sofic under approximations by 
supp-morphisms.  
\end{itemize} 
This proves that the set of triples of the statement of the theorem  
belongs to the class $\Pi_2$. 

To show that $\mathsf{Inf} = \{ m : W_m$ is infinite$\}$ is reducible 
to the family of triples $(D,(\circ ,e), n)$ as in the statement  
we define an algorithm which finds triples 
$(D_m , (\circ_m ,e_m ), n_m )$, where $(D_m ,\circ_m ) \cong (\mathbb{Z}/2\mathbb{Z} ,+)$, 
the number $n_m$ is a computability index of the function $x+3$ 
(i.e all $n_m$ are the same)  
and the members of $D_m$ are two elements $d_m$ , $e_m$ which 
satisfy the following property: 
\begin{quote} 
$e_m$ is a fixed index of the identity permutation of $\mathbb{N}$ 
and for all $m\in \mathsf{Inf}$ 
the number $d_m$ is a computability index of a permutation of order 2 with support $\mathbb{N} \setminus 3 \mathbb{N}$.  
\end{quote}  
At {\em Stage 0} we initialize the procedure by setting 
$f_{m,0}$ the function which for every natural $l$ takes  
the elements $3l, 3l+1$ to $3l+2$ and the elements $3l+2$ to $3l+1$. 
At {\em Stage} $s+1$ having input $(m, s+1)$ we verify if the algorithm 
enumerating $W_m$ adds a new element to the previously computed 
part of $W_m$, i.e. if $W_{m,s+1}\setminus W_{m,s}\not=\emptyset$. 
If this is not the case we set $f_{m,s+1} = f_{m,s}$. 
When $W_{m,s+1}\setminus W_{m,s}\not=\emptyset$ we 
change the definition of $f_{m,s}$ by $f_{m,s+1} (k) = k$, 
where $k$ is the first number of the form $3l$ which is taken by $f_{m,s}$ 
to $3l+2$. 

Let $f_m$ be the limit $\lim f_{m,s}$. 
It is clear that $f_m$ is computable and 
is the same for all $m$ with infinite $W_m$. 
Let $d_m$ be an index $t$ provided by the algorithm for $W_m$ such that 
$\varphi_t = f_m$.  
When $m\in \mathsf{Inf}$ the function $f_m$ is a permutation of order 2. 
Moreover in this case any restriction of 
$f_m$ to an initial segment of the form $\{ 0,1, \ldots , 3l\}$ 
is a permutation of order 2. 
Thus supp-morphisms of this form give sofic approximations. 
When $m\not\in \mathsf{Inf}$ the function $f_m$ is not a permutation, 
thus the triple $(D_m ,(\circ_m ,e_m ), n_m )$ 
does not satisfy the condition of the statement of the theorem.  
\end{proof}

\section{The computable part of $S_{\alpha}(\omega )$}

Let $\alpha$ be a growth which does not contain $\infty$ and which 
is represented by a computable function. 
Let 
\[
S^{comp}_{\alpha}(\omega ) = \{ \rho \in S(\omega ): \rho \mbox{ and } \rho^{-1} 
\mbox{ are computable and bounded above} 
\] 
\[
\mbox{ by some } g \in \alpha\}.
\] 
It is clear that 
\[
SF(\omega )< S^{comp}_{\alpha}(\omega )< S^{comp}(\omega ),   
\]
where $S^{comp}(\omega )$ consists of all computable permutations. 

Is it possible to get a version of Theorem \ref{3_2} 
where the numbering of computable functions $\varphi_k$ 
is replaced by a numbering $\nu$ of $S^{comp}_{\alpha}(\omega)$ 
which makes it a computably enumerable group?  
We will see below that $S^{comp}_{\alpha}(\omega)$ 
is not even computably enumerable.

\subsection{Computably enumerable actions} 

Let $G$  be a computably enumerable subgroup of $ S^{comp} (\omega)$ and 
let $\nu$ be the corresponding numbering.  
Here we use terminology of Section 2.1. 
Let $\mathsf{MultT}(G, \nu )$ be the graph of the multiplication. 
The natural action of $G$ on $\mathbb{N}$ defines the following ternary relation on $\mathbb{N}$: 
\[ 
R_{act}= \{ (l,m,n) : (\nu (l) )(m) = n \} . 
\]
When $R_{act}$ is computably enumerable it is computable: given $l,m,n$ 
one computes $(\nu(l))(m)$ and verifies $(\nu(l))(m)=n$. 
Thus assuming that $R_{act}$ is computably enumerable we arrive at the situation of a {\bf computable action}. 
The following proposition (Proposition \ref{pro1}) demonstrates that this condition makes some families of chunks slightly less complicated than 
in the case of the canonical numbering $\varphi_k$, $k\in \mathbb{N}$. 
Note that in Proposition \ref{pro1} the condition of statement 3 together with the condition of statement 2 quantified by $( \forall n )$ are exactly conditions provided by Theorem \ref{realization} for a sofic chunk. 
In fact Proposition \ref{pro1} describes their complexity in the case of computably enumerable groups. 
 
\begin{proposition} \label{pro1} 
Let $(G, \nu )$ be a computably enumerable subgroup of 
$S^{comp} (\omega )$ with the computable action 
and let $n_0$ be the number of the identity function on $\mathbb{N}$. 
Let $g: \mathbb{N} \rightarrow \mathbb{N}$ be a slow computable function 
from $\Omega_{\infty}$.  

\begin{enumerate} 
\item The following family is computable: \\ 
all triples $(D, \circ ,n_0 )$, where $D$ is a finite  subset of $\mathbb{N}$, $n_0 \in D$ and $\circ$ is a partial binary operation such that $E= \{ \nu (k) : k\in D\}$ is a chunk of $G$ which realizes 
$(D, \circ ,n_0 )$ with respect to the map $\nu (k) \rightarrow k$.   
\item The following family is computably enumerable: \\  
all quadruples $(D,\circ, n_0 , n)$, where the triple $(D, \circ ,n_0 )$ and the corresponding $E$ satisfy the conditions of statement 1 and there is $m\in \mathbb{N}$ and a $(1-\frac{1}{n} )$-expansive 
$\frac{1}{n}$-morphism $\sigma_m : E \rightarrow S_{m}$ such that 
\[ 
\forall k\in D \forall l\le m (\nu(k) (l) \le m \rightarrow  
\nu (k) (l) = \sigma_m (\nu (k) )(l))  .
\]  
\item  
The following family belongs to $\Pi_1$: \\ 
all triples $(D, \circ ,n_0 )$ as in statement 1, where 
$\{ \nu (k): k\in D\}$ is a $g$-chunk of $G$. 
\end{enumerate} 
\end{proposition} 

\begin{proof} 
Concerning statement 1 notice that computable enumerability of $(G, \nu )$
gives a computable enumeration of the family of all triples 
$(D, \circ ,n_0)$ where the map $\nu (k) \rightarrow k$ is not a homomorphism: given enumerations of $\mathsf{MultT}(G, \nu )$ and the family of all finite $(D, \circ )$  we form a list of such $(D, \circ ,n_0 )$ where some equalities from $\mathsf{MultT}(G, \nu )$ are not preserved. 
Therefore we see that the set of triples of the condition of statement 1 belongs to $\Pi_1$. 
On the other hand since the action of $(G, \nu )$  is computable there is a computable enumeration of the family of all triples 
$(D, \circ , n_0)$ where the map $\nu (k) \rightarrow k$ 
is a 1-1-homomorphism. 
See the beginning of the proof of Theorem \ref{3_2} for an analogous argument. 

Note that the set of all tuples $(D, \circ , n_0 , n,m )$ which satisfy the corresponding condition of statement 2 is computable. 
Thus applying $\exists m$ we obtain a computably enumerable set.  

To prove statement 3 of the proposition we use statement 1 and the assumption that $R_{act}$ is computable. 
The latter implies that the graph of each $\nu (k)$, $k\in D$, is 
computable. 
Thus stating that $\{ \nu (k) : k\in D\}$ is a $g$-chunk we may use only universal quantifiers. 
\end{proof}  

We will see in the next subsection that in this proposition  
$G$ cannot be $S^{comp}_{\alpha}(\omega)$. 
On the other hand for $G = S^{comp}_{\alpha}(\omega)$ 
relativized statements of this proposition still make sense:  
one can fix a Turing degree ${\bf d}$ and 
consider  computably enumerable groups and computably enumerable actions 
with respect to ${\bf d}$.

\subsection{The computable part of $S_{\alpha}(\omega )$}

The following theorem says that $S^{comp}_{\alpha}(\omega)$
is not computably enumerable. 

\begin{theorem} 
For any growth $\alpha$ the group 
$S^{comp}_{\alpha}(\omega )$ is not computably enumerable. 
\end{theorem} 

\begin{proof} 
This statement is known for the group $S^{comp}(\omega )$ (i.e. when $\alpha$ is represented by $\infty$), see p. 301, Section 6.1 in \cite{ershov}.  
We will adapt the corresponding proof. 
Since \cite{ershov} is not easily available and the proof given there needs small repairs we include some details. 

We firstly note that any finitary permutation belongs 
to $S^{comp}_{\alpha}(\omega)$ for any $\alpha$. 
We use the convention that for permutations $\sigma_1$ and $\sigma_2$ 
the permutation $\sigma_1\sigma_2$ acts by  $\sigma_2 (\sigma_1 (x))$. 
We will concentrate on transpositions. 
We need the following straightforward property of them:  
\begin{itemize} 
\item the supports of transpositions $\gamma$ and $\gamma'$ have a common 
point if and only if $(\gamma \cdot \gamma')^3 = 1$.    
\end{itemize} 
We will use the transpositions $\delta_1 = (0,1)$, $\delta_2 = (0,2)$. 
Let us consider the conjugacy class of this pair: 
$[ (\delta_1 , \delta_2 ) ]_{\alpha} = \{ (\delta^{\rho}_1 , \delta^{\rho}_2 ): \rho \in S^{comp}_{\alpha} (\omega)\}$. 
It is easy to see that it consists of all pairs of transpositions 
which supports have exactly one common point. 
Let us define the following relation on $[(\delta_1 , \delta_2 )]_{\alpha}$: 
\[ 
(\gamma_1 ,\gamma_2 ) \sim_{\eta} (\gamma'_1 ,\gamma'_2 )
\Leftrightarrow  \big(
(\gamma_1 \cdot \gamma'_1)^3 =  (\gamma_1 \cdot \gamma'_2)^3 = 
(\gamma_2 \cdot \gamma'_1)^3 = (\gamma_2 \cdot \gamma'_2)^3 = 1 \big)\wedge
\] 
\[
\big( \{ \gamma'_1 ,\gamma'_2 \} \not= \{ \gamma_1 , \gamma^{\gamma_1}_2 \} \big) \wedge \big( \{ \gamma'_1 ,\gamma'_2 \} \not= \{ \gamma_2 , \gamma^{\gamma_2}_1 \} \big) . 
\]
We claim that if pairs $(\gamma_1 ,\gamma_2 )$ and $(\gamma'_1 ,\gamma'_2 )$ 
belong to $[ (\delta_1 , \delta_2 ) ]_{\alpha}$, then they have the 
same intersections of supports if and only if 
$(\gamma_1 ,\gamma_2 ) \sim_{\eta} (\gamma'_1 ,\gamma'_2 )$. 
Indeed, the necessity is straightforward. 
For sufficiency assume $(\gamma_1 ,\gamma_2 ) \sim_{\eta} (\gamma'_1 ,\gamma'_2 )$. 
Then each pair of the family $\{ \gamma_1 ,\gamma_2 ,  \gamma'_1 ,\gamma'_2 \}$ has common points in supports. 
The inequalities of the definition of $(\gamma_1 ,\gamma_2 ) \sim_{\eta} (\gamma'_1 ,\gamma'_2 )$ say that when $\{ \gamma_1 ,\gamma_2 \} \not= \{ \gamma'_1 ,\gamma'_2 \}$, then the symmetric difference of 
$\mathsf{supp}(\gamma_1 )\cup \mathsf{supp} (\gamma_2 )$ and  $\mathsf{supp}(\gamma'_1 ) \cup \mathsf{supp}(\gamma'_2 )$ has at least two  elements. 
This condition forces the statement. 

The following permutation $\delta$ belongs to $S^{comp}_{\alpha}(\omega)$ 
for any $\alpha$: 
\[  
\delta (1) = 0 \mbox{ , } \delta (2n) = 2(n+1) \mbox{ for } n \in \mathbb{N} 
\mbox{ , and } \delta (2n +3 ) = 2n+1 \mbox{ for } n \in \mathbb{N}. 
\]  
Indeed both $\delta$ and $\delta^{-1}$ are bounded by $n+3$, 
which is the third iteration of $n+1$. 

For each $i$ let $(\gamma_{2i} , \gamma'_{2i})$ be 
$\delta^{-i} (\delta_1 , \delta_2 ) \delta^i$ and  
$(\gamma_{2i+1} , \gamma'_{2i+1})$ be 
$\delta^{i+1} (\delta_1 , \delta_2 ) \delta^{-i-1}$. 
Then $\mathsf{supp}(\gamma_{n}) \cap \mathsf{supp}(\gamma'_{n})=\{ n \}$. 
Moreover all $\mathsf{supp}(\gamma_i ) \cup \mathsf{supp}(\gamma'_i )$ 
are the following triples: 
\[ 
\ldots ,\{ 7,5,3\}, \{ 5,3,1\},\{ 3,1,0\},\{ 1,0,2\},\{ 0,2,4\},\{ 2,4,6\},\{ 4,6,8\}, \ldots . 
\] 
For for any permutation $\rho \in S^{comp}_{\alpha}(\omega)$ 
we have the following equivalence: 
\[  
\rho (k) = n \Leftrightarrow (\gamma_n , \gamma'_n )\sim_{\eta}  ( \gamma^{\rho}_k , (\gamma'_k)^{\rho}) . 
\] 
In particular: 
\begin{itemize} 
\item 
the equalities   
\[  
(\gamma_n \cdot \gamma^{\rho}_k)^3 = (\gamma'_n \cdot \gamma^{\rho}_k)^3 = 
(\gamma_n \cdot (\gamma'_k)^{\rho})^3 = (\gamma'_n \cdot (\gamma'_k)^{\rho})^3 = 1. 
\] 
imply that  
$\rho (k) \in \mathsf{supp} (\gamma_n ) \cup \mathsf{supp} (\gamma'_n ) $.
\end{itemize} 

Assuming that $S^{comp}_{\alpha}(\omega)$
is computably enumerable under a numbering $\nu$ we define an element of 
$S^{comp}_{\alpha}(\omega)$, say $\tau$, as follows. 
For every natural $n$ the permutation $\tau$ coincides with a 6-element cycle on $\{ 6n , 6n +1, 6n+2, 6n+3, 6n+4, 6n+5\}$. 
Thee are two cases. 
The first one is when $\rho = \nu (n)$ satisfies the following equalities:  
\[  
(\gamma_{6n +3}\cdot \gamma^{\rho}_{6n})^3 = (\gamma'_{6n +3} \cdot \gamma^{\rho}_{6n})^3 = 
(\gamma_{6n +3}\cdot (\gamma'_{6n})^{\rho})^3 = (\gamma'_{6n+3} \cdot (\gamma'_{6n})^{\rho})^3 = 1 .   
\] 
Then we put 
\[ 
\tau (x) = x+2 \mbox{ ( } \mathsf{mod} \mbox{ 6 ), where }
x \in \{ 6n, 6n+1 , 6n+2 , 6n+3 , 6n+4 , 6n+5 \} . 
\] 
In the contrary case we put $\tau (x) = x+3$ ( $\mathsf{mod}$ 6 ), 
where $x\in \{ 6n, 6n+1,  6n+2 , 6n+ 3 , 6n+4 , 6n+5\}$.
As a result $\tau$ belongs to $S_{\alpha}(\omega)$ and 
does not have any number under $\nu$. 
To see the latter note that when $\tau = \nu (n)$ then either 
$\nu (n)$ is as in the first case above but 
$\tau (6n ) = 6n + 2 \not\in \mathsf{supp}(\gamma_{6n +3})\cup \mathsf{supp} ( \gamma'_{6n+3} )$ or $\nu (n)$ is as in the second case but $\tau (6n) = 6n+3$, which is impossible. 

To get a contradiction with existence $\nu$  making 
$S^{comp}_{\alpha} (\omega )$ c.e. it remains to  notice that $\tau$ is a computable function.  
Note that since $(S^{comp}_{\alpha}(\omega), \nu)$ is computably enumerable 
there is an algorithm which lists all triples $k,m,n$ such that for 
$\rho = \nu (m)$ the following equalities are satisfied:  
\[  
(\gamma_n \cdot \gamma^{\rho}_k)^3 = (\gamma'_n \cdot \gamma^{\rho}_k)^3 = 
(\gamma_n \cdot (\gamma'_k)^{\rho})^3 = (\gamma'_n \cdot (\gamma'_k)^{\rho})^3 = 1. 
\] 
As a result it is decidable whether $\tau (6n) = 6n+2$. 
In particular $\tau$ is computable. 
\end{proof}  

\bigskip 

\begin{quote} 
{\bf Acknowledgements} 
The author is grateful to the anonymous referee for his thorough and essential remarks.   
\end{quote} 


\bigskip 

A. Ivanov  \\ 
Department of Applied Mathematics, Silesian University of Technology,\\  ul.Kaszubska 23, 44-101 Gliwice, Poland \\  
Aleksander.Iwanow@polsl.pl

\end{document}